\documentclass{article}
\usepackage{amssymb}
\usepackage{latexsym}
\newtheorem{thm}{Theorem}

\newtheorem{lem}{Lemma}
\newtheorem{prop}{Proposition}
\newtheorem{cor}{Corollary}

\begin{document}
\title{Holomorphic vector fields and minimal Lagrangian submanifolds}
\author{Edward Goldstein}
\maketitle

\renewcommand{\abstractname}{Abstract}
\begin{abstract}
The purpose of this note is to establish the following theorem: Let $N$ be a Kahler manifold, $L$ be an oriented immersed minimal Lagrangian submanifold of $N$ without boundary and $V$ be a holomorphic vector field in a neighbourhood of $L$ in $N$. Let $div(V)$ be the (complex) divergence of $V$. Then the integral $\int_{L}div(V) = 0$. Vice versa suppose that $N^{2n}$ is Kahler-Einstein with non-zero scalar curvature and $L$ is an embedded totally real $n$-dimensional oriented real-analytic submanifold of $N$ s.t. for any holomorphic vector field $V$ defined in a neighbourhood of $L$ in $N$, $\int_{L}div(V)=0$. Then $L$ is a minimal Lagrangian submanifold of $N$.
\end{abstract}
\section{Basic properties}
Let $(N^{2n},\omega)$ be a Kahler manifold. In this section we will discuss holomorphic vector fields on $N$ and present some basic facts on minimal Lagrangian submanifolds of $N$. The results of this section are essentially known, though they might be stated in different terms in the literature.

First we discuss holomorphic vector fields on $N$. Let $V$ be  a vector field defined on some open subset $U$ of $N$. The following proposition is elementary and well known (see \cite{Kob}, Proposition 4.1):
\begin{prop}
The following conditions are equivalent:

1) The flow of $V$ commutes with the complex structure $J$ on $N$.

2) For any point $m \in U$ the following endomorphism : $X \mapsto \nabla_X V$ of $T_mN$ is $J$-linear on $T_m N$.

3) The vector field $V-iJV$ gives a holomorphic section of $T^{(1,0)}U$.
\end{prop}
A vector field $V$ satisfying the conditions of Proposition 1 is called a {\it holomorphic} vector field. Let $V$ be a holomorphic vector field on some open subset $U$ of $N$ and let $m$ be a point in $U$. Since the endomorphism $X \mapsto \nabla_X V$ is $J$-linear on $T_mN$ we can define 
\begin{equation}
div(V)=trace_{\mathbb{C}}(X \mapsto \nabla_X V)
\end{equation}
Let $f=Ref + i Imf$ be a holomorphic function on $U$. From condition 3) of Proposition 1 we deduce that the vector field $fV=Ref V + ImfJV$ is a holomorphic vector field on $U$. Moreover one easily computes that 
\begin{equation}
div(fV)=fdiv(V)+V(f)
\end{equation}
Let $K(N)$ be the canonical bundle of $N$ (i.e. $K(N)=\Lambda^{(n,0)}T^{\ast}N$). Let $\varphi$ be a section of $K(N)$ over $U$ (not necessarily holomorphic). Thus $\varphi$ is an $(n,0)$-form over $U$. 
\begin{prop}
Let $V$ be a holomorphic vector field on $U$. Then 
\[ \nabla_V \varphi= {\cal L}_V \varphi- div(V)\varphi \] 
\end{prop}
{\bf Proof:} Let $m \in U$. Pick a unitary basis $X_1,\ldots,X_n$ of $T_mN$ (here $T_mN$ is viewed as Hermitian vector space with the complex structure $J$). Extend $X_i$ to a unitary frame in a neighbourhood of $m$. Then
\begin{equation}
{\cal L}_V \varphi(X_1,\ldots,X_n)=V(\varphi(X_1,\ldots,X_n))- \Sigma \varphi(X_1,\ldots,[V,X_n],\ldots,X_n)
\end{equation}
and 
\begin{equation}
\nabla_V \varphi(X_1,\ldots,X_n)=V(\varphi(X_1,\ldots,X_n))-\Sigma \varphi (X_1,\ldots,\nabla_V X_i,\ldots,X_n)
\end{equation} 
Now $\nabla_V X_i= [V,X_i] + \nabla_{X_i}V$. We plug this into (4) and subtract (4) from (3) to deduce the statement of the proposition. Q.E.D.

Next we prove the following lemma (which we essentially proved in \cite{Gold}):
\begin{lem}
Let $(N,\omega)$ be a Kahler-Einstein manifold with non-zero scalar curvature $t$. Let $V$ be a holomorphic infinitesimal isometry on some neighbourhood $U$ of $N$. Then the function $\mu=it^{-1}div(V)$ is a moment map for the $V$-action on $(N,\omega)$
\end{lem}
{\bf Proof:} We need to prove that $d\mu=i_V \omega$. We shall prove it at a point $m$ s.t. $V(m) \neq 0$. Pick an element $\varphi$ of $K(N)$ over $m$ which has unit length. Since the flow of $V$ is given by holomorphic isometries we can extent $\varphi$ to a unit length section of $K(N)$ invariant under the $V$-flow on some neighbourhood $U$ of $m$. The section $\varphi$ defines a connection 1-form $\xi$ on $U$, $\xi(u)=<\nabla_u\varphi,\varphi>$. The Einstein condition tells that 
\begin{equation} 
id\xi=t \omega 
\end{equation}
Since $\varphi$ is $V$-invariant we deduce from Proposition 2 that 
\begin{equation}
div(V)=-\xi(V)
\end{equation}
Also since $\varphi$ is $V$-invariant and the flow of $V$ is given by isometries, we deduce that $\xi$ is also $V$-invariant. Thus \[0={\cal L}_V \xi= d(\xi(V))+ i_V d\xi=~ by ~ (5) ~ and ~ (6)= -d(div(V)) -iti_V\omega \]
and the lemma follows. Q.E.D.

Next we discuss minimal Lagrangian submanifolds on $N$. Let $L$ be an oriented $n$-dimensional totally real submanifold of $N$ (i.e. $TL \bigcap J(TL)=0$). For any point $l \in L$ there is a unique element $\kappa_l$ of $K(N)$ over $l$ which restricts to the volume form on $L$. Various $k_l$ give rise to a section 
\begin{equation}
\kappa:L \mapsto K(N)
\end{equation}
Let now $L$ be a Lagrangian submanifold of $N$. The section $\kappa$ is a unit length section of $K(N)$ over $L$ and it
defines a connection 1-form $\xi$ for the connection on $K(N)$ over $L$,
$\xi(u)= <\nabla_u \kappa, \kappa>$. Here $\nabla$ is the connection on
$K(N)$, induced from the Levi-Civita connection on $N$. Since $\kappa$ has unit length $\xi$ is an imaginary valued 1-form.
 
Let $h$ be the
trace of the second fundamental form of $L$. So $h$ is a section of the
normal bundle of $L$ in $N$ and we have a corresponding 1-form $\sigma=
i_h\omega$ on $L$. The following fact is well-known, although it is often
stated differently in the literature (see \cite{Sch}):
\begin{lem} 
$\sigma= i \xi$
\end{lem}
{\bf Proof:} Let $l \in L$ and $e$ be some vector in the tangent space to
$L$ at $L$. To compute $\xi(e)$ we need to compute $\nabla_e \kappa$. Take
an orthonormal frame $(v_j)$ of $T_lL$ and extend it to an orthonormal frame in a neighbourhood $U$ of $l$ in $L$ s.t. $\nabla^L v_i=0$ at $l$ (here $\nabla^L$ is the Levi-Civita connection of $L$). We get that \[ \nabla_e \kappa = \kappa \cdot \nabla_e \kappa(v_1,\ldots,v_n)= \kappa (e(\kappa(v_1,\ldots,v_n))- \Sigma\kappa(v_1,\ldots,\nabla_e v_j,\ldots,v_n)) \]
Now $e(\kappa(v_1,\ldots,v_n))=0$. Also clearly \[\kappa(v_1,\ldots,\nabla_e v_j,\ldots,v_n)=i<\nabla_e v_j,Jv_j>=i<\nabla_{v_j}e,Jv_j>=i<-e,J(\nabla_{v_j}v_j)> \]   
Here $J$ is the complex structure on $N$. Thus we get that 
\[\nabla_e \kappa = -i(Jh \cdot  e)\kappa_l = -i \sigma(e)\kappa_l\] 
Here $h= \Sigma \nabla_{v_j}v_j$ is the trace of the second fundamental
form of $L$. Thus $\sigma= i \xi$.  Q.E.D.

Thus if $L$ is {\it minimal} (i.e. $h=0$) iff $\kappa$ is parallel over $L$.

{\bf Remark:} Let $L$ be a minimal Lagrangian submanifold of $N$. We have seen that $\xi=0$ on $L$. Thus also $d\xi=0$ on $L$. But $d\xi$ is the curvature form for the connection on $K(N)$, i.e. $d\xi=-iRic$. Here $Ric$ is the Ricci form of $N$, and $Ric$ is proportional to $\omega$ iff $N$ is Kahler-Einstein. If $N$ is Kahler-Einstein the condition $Ric|_L=0$ follows from the Lagrangian condition on $L$. But if $N$ is not Kahler-Einstein, we have a new algebraic condition $Ric|_L=0$ on minimal Lagrangian submanifolds.
\section{Proof of the main theorem}
We now can state and prove our main theorem:
\begin{thm}

1) Let $N$ be a Kahler manifold, $L$ be an oriented immersed minimal Lagrangian submanifold of $N$ without boundary and $V$ be a holomorphic vector field defined in a neighbourhood of $L$ in $N$. Then \[\int_{L}div(V)=0 \]

2) Let $N^{2n}$ be a Kahler-Einstein manifold with non-zero scalar curvature and $L$ be an $n$-dimensional totally real oriented embedded real-analytic submanifold of $N$ s.t. for any holomorphic vector field $V$ defined in a neighbourhood of $L$ in $N$ we have $\int_{L}div(V)=0$. Then $L$ is a minimal Lagrangian submanifold of $N$
\end{thm}
{\bf Proof:} 1) Let $L$ be a minimal Lagrangian submanifold of $N$ and $V$ be a holomorphic vector field defined in a neighbourhood of $L$ in $N$. Let $\kappa$ be a section of $K(N)$ over $L$ as in equation (7). Since $\kappa$ restricts to the volume form on $L$ we have $\int_{L}div(V)=\int_{L}div(V)\kappa$. Let 
\begin{equation}
\phi=i_V\kappa|_L
\end{equation}
$\phi$ is an $(n-1)$-form on $L$. We claim that 
\begin{equation}
d\phi=div(V)\kappa|_L
\end{equation}
Thus the first assertion of the theorem will follow. To prove (9) let $l$ be a point in $L$. By Lemma 2 we have that for any element $w$ in the tangent bundle to $L$, $\nabla_w \kappa=0$. We can extend $\kappa$ to be a section of $K(N)$ over some neighbourhood $Z$ of $l$ in $N$ s.t. for any element $w$ in the normal bundle of $L$ to $N$ in $Z \bigcap L$ we'll have $\nabla_w \kappa=0$. Thus we'll have $\nabla \kappa= 0$ along $L$. From this it also follows that $d\kappa=0$ along $L$.
Now we use equation Proposition 2 for $V$ and $\varphi=\kappa$. We deduce that 
\[ div(V)\kappa= {\cal L}_V \kappa \]
along $L$. Also ${\cal L}_V \kappa=d(i_V \kappa)+i_V(d\kappa)$ and $d\kappa$ vanishes along $L$. Thus we get \[div(V)\kappa|_L=d \phi \]

2) Let $(N^{2n},\omega)$ be a Kahler-Einstein manifold with a non-zero scalar curvature $t$ and $L$ be an oriented real-analytic embedded totally real $n$-dimensional submanifold of $N$ s.t. $\int_{L}div(V)=0$ for any holomorphic vector field $V$ near $L$. Consider the section $\kappa$ of $K(N)$ over $L$ as in (7). Since $L$ is totally real, $n$-dimensional embedded real-analytic submanifold of $N$, one can uniquely extend $\kappa$ to a holomorphic section of $K(N)$ over some neighbourhood $U'$ of $L$ in $N$ (see the Appendix). 

Let $\xi=Re\xi+iIm \xi$ be the connection 1-form on $L$ defined by the section $\kappa$ over $L$, i.e. for any tangent vector $u$ to $L$ we have $\nabla_u\kappa=\xi(u)\kappa$. Let $V_r$ be the vector field on $L$ dual to the form $Re \xi$ with respect to the Riemannian metric on $L$. $V_r$ is a real-analytic vector field on $L$ and by Proposition 4 of the Appendix we can extend $V_r$ to a holomorphic vector field $V_r$ on a neighbourhood of $L$ in $N$. By Proposition 2 we have 
\[ \xi(V_r)\kappa= {\cal L}_{V_r}\kappa -div(V_r)\kappa ~ on ~ L\]
We integrate this over $L$ to get
\begin{equation}
\int_{L}\xi(V_r)\kappa=\int_{L}{\cal L}_{V_r}\kappa - \int_{L}div(V_r)\kappa
\end{equation}  
Now ${\cal L}_{V_r}\kappa=d(i_{V_r}\kappa)$ and it integrates to $0$ over $L$. Also by our assumptions since $\kappa$ restricts to the volume form on $L$ we get \[\int_{L}div(V_r)\kappa= \int_{L}div(V_r)=0\]
Thus $\int_{L}\xi(V_r)\kappa=0$. But $Re(\xi(V_r))=|V_r|^2$ pointwise. Thus $V_r=0$ and so $Re\xi=0$. Similarly we prove that $Im \xi=0$.

Thus $\xi=0$ on $L$. So $d\xi=0$ on $L$. But \[d\xi=-it \omega|_L \] Here $\omega$ is the Kahler form on $N$. Hence $L$ is Lagrangian. Since $\xi=0$ on $L$ we deduce from Lemma 2 that $L$ is minimal. Q.E.D.

Let us derive a simple corollary of Theorem 1:
\begin{cor}
Let $L$ be an immersed oriented minimal Lagrangian submanifold of $\mathbb{C}P^n$ and let $(z_1, \ldots, z_{n+1})$ be homogeneous coordinates on $\mathbb{C}P^n$. Then we can't have $|z_1| > |z_2|$ at all points of $L$.
\end{cor}
{\bf Proof:} Consider the following circle action on $\mathbb{C}P^n$:
\[ e^{i\theta}(z_1, \ldots,z_{n+1})=(e^{i \theta}z_1,e^{-i \theta}z_2,z_3, \ldots, z_{n+1}) \]
Let $V$ be the vector field on $\mathbb{C}P^n$ generating this action. 
$\mathbb{C}P^n$ is Kahler-Einstein with scalar curvature 1, hence by Lemma 1 the function $idiv(V)$ is a moment map for the $S^1$-action on $\mathbb{C}P^n$.
We have computed in \cite{Gold} that 
\[idiv(V)= (|z_1|^2-|z_2|^2)/\Sigma|z_i|^2 \]
In fact we can also deduce this from Theorem 1.
Indeed the map $f=(|z_1|^2-|z_2|^2)/\Sigma|z_i|^2$ is a moment map for the $S^1$-action on $\mathbb{C}P^n$, hence it differs from $idiv(V)$ by a constant $c$. Also the submanifold $L'=((z_1,\ldots,z_{n+1})||z_1|=|z_j|)$ is a minimal Lagrangian submanifold of $\mathbb{C}P^n$. Hence by Theorem 1 $\int_{L'}div(V)=0$. From this we deduce that $c=0$ i.e. $idiv(V)=f$.

Let now $L$ be an immersed oriented minimal Lagrangian submanifold of $\mathbb{C}P^n$.
We have $\int_{L}div(V)=0$. Hence we obviously can't have $|z_1|>|z_2|$ everywhere on $L$. Q.E.D.  
\section{Appendix}
In this Appendix we want to demonstrate the following fact (used in the proof of Theorem 1): Let $L$ be a totally real $n$-dimensional embedded real-analytic compact submanifold of a complex manifold $N^{2n}$. Suppose that $P$ is a holomorphic vector bundle over $N$ and we have a real-analytic section $\sigma$ of $P$ over $L$. Then we can uniquely extend it to a holomorphic section $\sigma'$ over some neighbourhood of $L$ in $N$. We begin with the following proposition:
\begin{prop}
Let $f : U \mapsto \mathbb{C}^k$ be a real analytic map from an open subset $U$ of $0$ in $\mathbb{R}^n$ to $\mathbb{C}^k$. Then we can uniquely extend $f$ to a holomorphic map $f'$ from some open subset $U'$ of $0$ in $\mathbb{C}^n$ to $\mathbb{C}^k$. Here we think of $\mathbb{R}^n$ as a subset of $\mathbb{C}^n$.
\end{prop}
{\bf Proof:} Let $(z_1, \ldots, z_k)$ be coordinates on $\mathbb{C}^k$. We can think of $f$ as \[f=(f_1,\ldots,f_k) \]
and we need to extend each $f_i$ to a holomorphic function on an open subset of $0$ in $\mathbb{C}^n$. Let $x=(x_1,\ldots,x_n)$ be the coordinates on $\mathbb{R}^n$. Since $f_i$ is real-analytic on $\mathbb{R}^n$ we can write its Taylor's expansion \[ f_i= \Sigma C_{\alpha}x^{\alpha} \]
near $0 \in \mathbb{R}^n$. Clearly $f_i$ has a unique holomorphic extension \[f'_i=\Sigma C_{\alpha}z^{\alpha}\]
onto a neighbourhood of $0$ in $\mathbb{C}^n$. Q.E.D.

Now we can prove the main result of the Appendix:
\begin{prop}
Let $L$ be a totally-real $n$-dimensional embedded compact real-analytic submanifold of a complex manifold $N$. Suppose that $P$ is a holomorphic vector bundle over $N$ and we have a real-analytic section $\sigma$ of $P$ over $L$. Then $\sigma$ extends uniquely to a holomorphic section $\sigma'$ on a neighbourhood of $L$ is $N$
\end{prop}
{\bf Proof:} It is obviously enough to prove that for any point $l \in L$ we can uniquely extend $\sigma$ onto a neighbourhood of $l$ in $N$. Also near $l$ we can think of $P$ as being the trivial bundle $\mathbb{C}^k$. Suppose now that there is a biholomorphic map $\phi$ from a neighbourhood $U$ of $0$ in $\mathbb{C}^n$ onto a neighbourhood $U'$ of $l$ in $N$ s.t. $\phi(\mathbb{R}^n \bigcap U) = L \bigcap U'$. Then the desired claim will follow from Proposition 3.

To construct the biholomorphic map $\phi$ we again use Proposition 3. Since $L$ is a real-analytic submanifold of $N$ we can find a neighbourhood $W$ of $0$ in $\mathbb{R}^n$, a neighbourhood $W'$ of $l$ in $L$ and a real-analytic map $f:W \mapsto N$ s.t. the image of $f$ lies in $L$ and in fact $f : W \mapsto W'$ is a diffeomorphism. Since $N$ is a complex manifold we can find a neighbourhood of $l$ in $N$, which is biholomorphic to a ball in $\mathbb{C}^n$. Thus we can think of $f$ as a map form $W$ to $\mathbb{C}^n$. By Proposition 3 we can extend it to a holomorphic map $f'=\phi$ from a neighbourhood of $0$ in $\mathbb{C}^n$ to $\mathbb{C}^n$. Since $f:W \mapsto W'$ is a diffeomorphism it is clear that the differential of $\phi$ at $0 \in \mathbb{C}^n$ is an isomorphism. Thus $\phi$ is a biholomorphic map from some neighbourhood $U$ of $0$ in $\mathbb{C}^n$ and we are done. Q.E.D.

Massachusetts Institute of Technology

E-Mail : egold@math.mit.edu


\begin {thebibliography}{99}  
\bibitem[1]{Gold} E. Goldstein : Calibrated Fibrations on complete
manifolds via torus action, math.DG/0002097
\bibitem[2]{Kob} S. Kobayashi : Transformation groups in differential geometry, Ergebnisse der Mathematik und ihrer Grenzgebiete, Bd. 70, Springer-Verlag 1972
\bibitem[3]{Sch} R. Schoen and J. Wolfson: Minimizing volume among Lagrangian submanifolds, Differential Equations: La Pietra 1996 (Florence), 181-199
\end{thebibliography}
\end{document}